\documentclass[12pt]{article}
\usepackage{amsfonts,amsmath,amssymb,mathrsfs}
\title{Some Jump Processes in Quantum Field Theory}
\author{
Roderich Tumulka\footnote{Dipartimento di
        Fisica and INFN sezione di Genova, Universit\`a di Genova, Via
        Dodecaneso~33, 16146~Genova, Italy.
        \texttt{tumulka@mathematik.uni-muenchen.de}}\ \ and
Hans-Otto Georgii\footnote{Mathematisches Institut,
        Ludwig-Maximilians-Universit\"at, Theresienstr.~39,
        80333~M\"unchen, Germany.
        \texttt{georgii@mathematik.uni-muenchen.de}}}
\date{}

\newcommand{\E}{\mathrm{e}} 
\newcommand{\I}{\mathrm{i}} 
\newcommand{\D}{\mathrm{d}} 
\newtheorem{theorem}{Theorem}

\newcommand{\CCC}{\mathbb{C}} 
\newcommand{\RRR}{\mathbb{R}} 
\newcommand{\ZZZ}{\mathbb{Z}} 
\newcommand{\EEE}{\mathbb{E}} 
\newcommand{\1}{\mathbf{1}} 
\newcommand{\tr}{\mathrm{tr}} 
\newcommand{\Laplace}{\Delta} 
\newcommand{\im}{\mathrm{Im}} 

\renewcommand{\limits}{}

\newcommand{\vr}{x}
\newcommand{\vy}{y}

\newcommand{\vQ}{Q}

\newcommand{\valpha}{\alpha}

\newcommand{\Gommo}{\Gamma_{\!\neq}} 
\newcommand{\Hilbert}{\mathscr{H}} 
\renewcommand{\sp}[2]{\langle #1 | #2 \rangle} 
\newcommand{\conf}{\mathcal{Q}} 
\newcommand{\interior}{\conf^\circ} 
\newcommand{\boundary}{\partial \conf} 
\newcommand{\salg}{\mathscr{F}} 
\newcommand{\prob}{\mathbb{P}} 
\newcommand{\pov}{{P}} 
\newcommand{\domain}{\mathrm{domain}} %
\newcommand{\measure}{\pi} 
\newcommand{\current}{J} 
\newcommand{\generator}{\mathscr{L}} 
\newcommand{\friedhof}{\infty} 
\newcommand{\dest}{f} 
\newcommand{\area}{\lambda} 
\newcommand{\vol}{\mu} 
\newcommand{\antidest}{\nu} 

\begin{document}\maketitle
\begin{abstract}
      A jump process for the positions of interacting quantum particles
      on a lattice, with time-dependent transition rates governed by the
      state vector, was first considered by J.S.~Bell. We review this
      process and its continuum variants involving ``minimal'' jump
      rates, describing particles as they get created, move, and get
      annihilated. In particular, we sketch a recent proof of global
      existence of Bell's process.  As an outlook, we suggest how
      methods of this proof could be applied to similar global existence
      questions, and underline the particular usefulness of minimal jump
      rates on manifolds with boundaries.

\medskip

\end{abstract}

\section{Introduction}

This contribution deals with Markov jump processes $Q_t$ describing
the positional time evolution of finitely many interacting quantum
particles.  These processes are characterized by a specific form of
time-dependent jump rates induced by the Schr\"odinger equation for
the quantum state vector $\Psi_t$ of the underlying quantum field
theory (QFT).  As a typical example, suppose that the particles live
in the physical three-space $\RRR^3$. The process $Q_t$ then takes
values in the space of all finite subsets of $\RRR^3$, the
configuration space of a variable number of identical particles
(corresponding to the Fock space of QFT), and the jumps of $Q_t$
describe the creation or annihilation of particles; between these
stochastic jumps, $Q_t$ evolves deterministically according to some
ordinary differential equation governed by $\Psi_t$.  Alternatively,
one may think of quantum particles on a lattice; the jumps of $Q_t$
then record all changes of the particle configuration. We will portray
several processes of this type, present some common principles, and in
particular discuss some results of two recent papers
\cite{crex1,crea2}, the work on which was supported by the DFG
Priority Program.

As the state vector $\Psi_t$ determining the jump rates follows the
time-dependent Schr\"o\-din\-ger equation, the jump rates themselves
are explicitly time-dependent, so that the processes $Q_t$ considered
here do not admit an invariant measure.  However, the jump rates are
designed in such a way that $Q_t$ does admit an \emph{equivariant
measure}, namely the quantum distribution $|\Psi_t|^2$, which means
that $Q_t$ has distribution $|\Psi_t|^2$ at any time $t$. This is the
key feature justifying the particular form of the jump rates, and on
the other hand the main fact on which one can build the analysis of
these processes.  So, the issue here is not the analysis of
distributional properties of a given process, but the converse: the
equivariant distribution is given, and the objective is to prove the
existence of the associated process, and to check that it really does
have the equivariant distribution.  In \cite{crex1}, we carried out
this program for the case of a discrete configuration space, including
in particular a lattice model proposed by J.S.~Bell \cite{BellBeables};
the main arguments will be sketched in Sect.~4.

{}From the probabilistic viewpoint, one has to overcome two
difficulties.  First, the transition rates exhibit singularities, in
that they become ill-defined at certain time and space points. One has
to show that the process avoids these singularities. The second (and
more important) task is to rule out the possibility of explosion,
i.e., the accumulation of infinitely many jumps in finite time. Due to
the {unbounded growth of the rates near the} singularities, the
standard methods fail, and one has to use the particular relation
between transition rates and equivariant distribution.

Besides the results on the discrete case mentioned above, we will also
describe the continuum analogues of Bell's process investigated in
\cite{crea1,crea2}; as a special case these include Bohmian mechanics
\cite{Bohm52,Bellhidden,survey,DetlefBuch} for the appropriate
Hamiltonian with a conserved number of particles. On the basis of what
we learned from our existence proof for Bell's model, we also propose
here some new methods for proving the existence of Bohmian mechanics.

Let us now discuss how the models considered here relate to the topics
of other articles in this volume. First, the existence problem for a
model of quantum field theory is also the subject of the contribution
of S.~Albeverio, Y.~Kondratiev, M.~R\"ockner, and T.~Pasurek.
The issue there is the existence, and uniqueness, of Euclidean Gibbs
measures for infinitely many interacting quantum spins. These concern
an equilibrium setting, and time appears only via path integration to
make the connection with the quantum states. The difficulty there is
the infinite number of spins, requiring particular assumptions on the
interaction. By way of contrast, the models considered here involve
only finitely many particles, but in a nonequilibrium situation, and
we do not need any particular assumptions on the interaction.

{}From the methodological side, there is a closer connection to the
contributions dealing with population biology, in particular that by
R.~H\"opfner and E.~L\"ocherbach.  The similarities concern the
creation and annihilation of particles in the vicinity of other
particles, and the necessity of proving non-explosion.  Also, the
space-dependence of the reproduction rates in H\"opfner and
L\"ocherbach's article implies non-exponential life-times of
individuals, just as the time-dependent jump rates considered here
imply non-exponential interjump times.  However, in our case the paths
between the jumps are smooth and deterministic.

This note is organized as follows.  In Sect.~2 we derive the
fundamental formula \eqref{tranrates} for minimal jump rates, defining
the jump process associated with a certain type of Hamiltonian.  This
involves consideration of the equivariant probability distribution
\eqref{measuredef} and probability current \eqref{currentdef} provided
by quantum theory. In Sect.~3 we explain the connection with Bohmian
mechanics and with Bell's model. We also describe the process for a
concrete example QFT, introduced in \cite{crea1}. In Sect.~4 we sketch
the global existence proof for the discrete case, including Bell's
model, that we developed in \cite{crex1}. In Sect.~5 we point out how
the methods of \cite{crex1} could be adapted to other global existence
problems. In Sect.~6 we indicate some perspectives for future research
concerning a process for quantum theory on a manifold with boundaries,
and the special role the minimal jump rate \eqref{tranrates} plays for
this process.

\section{Jump Rates Induced by Schr\"odinger Equations}

We now introduce the class of jump processes we are concerned with,
starting with a general framework.  For our purposes, a quantum theory
is abstractly given by a Hilbert space $\Hilbert$ containing the state
vectors, a one-parameter group $U_t$ of unitary operators on
$\Hilbert$ defining the time evolution
\begin{equation}\label{evolution}
      \Psi_t = U_t \Psi_0
\end{equation}
of the state vector, and a measurable space $(\conf,\salg)$ of
configurations describing the locations of particles. $\conf$ is tied
to $\Hilbert$ via a projection-valued measure (PVM) $\pov(\D q)$ on
$\conf$ acting on $\Hilbert$, i.e., a mapping from the
$\sigma$-algebra $\salg$ to the family of projection operators on
$\Hilbert$ that is, like a measure, countably additive (in the sense
of the weak operator topology) and normalized, in that $\pov(\conf) =
I$, the identity operator on $\Hilbert$.  If $\Hilbert =
L^2(\conf,\CCC^k)$ with respect to some measure on $\conf$ then
$\Hilbert$ is equipped with a natural PVM, namely $\pov(B)$ being
multiplication by the indicator function of the set $B$. In
nonrelativistic quantum mechanics, another way of saying this is that
$\pov$ is the PVM corresponding to the joint spectral decomposition of
all position operators.

By Stone's theorem, the unitary operators $U_t$ are of the form
\begin{equation}\label{Hamiltonian}
      U_t = \E^{-\I Ht/\hbar}
\end{equation}
with $H$ a self-adjoint operator on $\Hilbert$, called the
Hamiltonian. Equations \eqref{evolution} and \eqref{Hamiltonian}
together
correspond to the formal Schr\"odinger equation
\begin{equation}\label{schr}
      \I\hbar\, \frac{\D \Psi_t}{\D t} = H \Psi_t\,.
\end{equation}
We will show how this Schr\"odinger equation, together with the PVM
$\pov$, gives rise to a natural Markov process on $\conf$.

In this section we focus on the case in which this Markov process is a
pure jump process. (Roughly speaking, this will require that the
Hamiltonian is an integral operator in the ``position representation''
defined by $\pov$; differential operators will be considered in
Sect.~3.) So we ask: \emph{Is there any distinguished Markovian jump
process $(Q_t)$ on $\conf$ describing the evolution of the particle
configuration, and what are its transition rates?}  To answer this
question we note that the quantum theoretical probability distribution
of the configuration at time $t$ is given by
\begin{equation}\label{measuredef}
      \measure_t(\,\cdot\,) = \sp{\Psi_t}{\pov(\,\cdot\,) \Psi_t}.
\end{equation}
(We generally assume that $\|\Psi_0\|=1$.) It is therefore natural to
stipulate that $\measure_t$ is \emph{equivariant} for $(Q_t)$, meaning
that $Q_t$ has distribution $\measure_t$ at every time $t$. Can one
choose some (time-dependent) transition rates $(\sigma_t)$ for $(Q_t)$
to satisfy this requirement of equivariance? Yes indeed, in view of
\eqref{schr} the time derivative of $\measure_t$ is given by
\begin{equation}\label{dotmeasure}
      \dot{\measure}_t(B) = \tfrac{2}{\hbar} \, \im \,
      \sp{\Psi_t}{ \pov(B) H \Psi_t}
      = \int  \current_t(B, \D q')\,,
\end{equation}
where
\begin{equation}\label{currentdef}
      \current_t(B, B') = \tfrac{2}{\hbar} \, \im \, \sp{\Psi_t}{
      \pov(B) H \pov(B') \Psi_t}
\end{equation}
is the quantum theoretical current between two sets $B,B'\in\salg$.
On the other hand, suppose $(Q_t)$ is a pure jump process on $\conf$
jumping at time $t$ with rate $\sigma_t(B,q')$ from $q'\in\conf$ to
some configuration in $B\in\salg$. Then its distribution
$\rho_t=\prob\circ Q_t^{-1}$ evolves according to the equation
\begin{equation}\label{dotrho}
     \dot{\rho}_t(B) = \int_{\conf}\sigma_t(B,q)\,
     \rho_t(\D q) - \int_B \sigma_t(\conf, q)\, \rho_t(\D q) .
\end{equation}
To satisfy the condition of equivariance we need to find jump rates
$\sigma_t$ such that the right-hand sides of the evolution equations
\eqref{dotmeasure} and \eqref{dotrho} coincide whenever
$\rho_t=\measure_t$. We see that this is the case when $\sigma_t$ is
given by the Radon--Nikodym derivative
\begin{equation}\label{tranrates}
      \sigma_t(\D q,q') =
      \frac{\current_t^+(\D q,\D q')}{\measure_t(\D q')}
      =\frac{\bigl[ \tfrac{2}{\hbar} \im \,
      \sp{\Psi_t}{\pov(\D q) H \pov(\D q') \Psi_t} \bigr]^+}
      {\sp{\Psi_t}{\pov(\D q') \Psi_t}}
\end{equation}
of the positive part $\current_t^+$ of $\current_t$ in its second
variable $q'$, provided this makes sense. For, the antisymmetry of
$\current_t$ then implies that
\[
   \sigma_t(\D q,q')\measure_t(\D q') - \sigma_t(\D q',q)\measure_t(\D q)
   =  \current_t(\D q, \D q')\,.
\]

To make formula \eqref{tranrates} meaningful one needs some
assumptions which roughly require that $H$ is an integral operator in
the position representation given by $\pov$, and $(\conf,\salg)$ is
standard Borel. This is discussed in detail in \cite{crea2}, where
\eqref{tranrates} was written down for the first time in this
generality; special cases had been used before in
\cite{BellBeables,Sudbery,crea1}. For the precise conditions we refer
to Theorem~1 (Sect.~4.1) and Corollaries~1--3 (Sect.~4.2) of
\cite{crea2}. Under these conditions, formula \eqref{tranrates} can
(and has to) be read as follows: A priori, $\current_t$ is a signed
bi-measure on $\salg \times\salg$ (a measure in each of the two
variables $q,q'$). This has to (and then can) be extended to a signed
measure on the product $\sigma$-algebra $\salg\otimes\salg$. The
positive part in \eqref{tranrates} is then to be understood in the
sense of the Hahn--Jordan decomposition of this extended measure. Next
one notes that, for each $B\in\salg$, $\current_t(B,\cdot) \ll
\measure_t$ because $\pov(B')\Psi_t=0$ whenever
$\measure_t(B')=0$. One can thus form the Radon--Nikodym derivatives
$\sigma_t(B,\cdot)=J_t^+(B,\D q')/\measure_t(\D q')$, which finally
have to be chosen in such a way that $\sigma_t$ becomes a measure
kernel.

According to our derivation above, the transition rates
\eqref{tranrates} have been chosen to satisfy the requirement of
equivariance. There was, however, still some freedom of choice. The
particular rate \eqref{tranrates} is singled out by the following
additional facts.

\begin{enumerate}
\item Suppose there exists a jump process $(Q_t)$ on
$\conf$ with rates \eqref{tranrates}. As is evident from the arguments
above,
the net probability current of $(Q_t)$ between two sets
$B,B'\in\salg$,
\[
      j_t(B, B') = \lim_{\varepsilon \searrow 0} \frac{1}{\varepsilon}
      \Bigl( \prob(Q_t \in B', Q_{t+\varepsilon} \in B) - \prob(Q_t \in
      B, Q_{t+\varepsilon} \in B') \Bigr),
\]
then coincides with the quantum theoretical current $\current_t$
defined by \eqref{currentdef}.  Conversely, if $(\tilde Q_t)$ is any
pure jump process having initial distribution $\measure_0$ at time
$0$, some jump rates $\tilde \sigma_t$ and probability current $\tilde
\jmath_t =\current_t$, it turns out that necessarily $\tilde \sigma_t
(\D q ,q') \geq \sigma_t (\D q ,q')$ \cite{Roy,BD,crea2}. This
follows from the minimality of the Hahn--Jordan decomposition.  The
rates \eqref{tranrates} are therefore called the \emph{minimal jump
rates}, and a process with rates \eqref{tranrates} is distinguished
among all processes with the right probability current by having the
least frequent jumps, or the smallest amount of randomness.

\item Always one of the transitions $q' \to q$ or $q\to q'$ is
forbidden. More precisely, for every time $t$ there exists a set
$S_t^-\in\salg\otimes \salg$ which, together with its transposition
$S_t^+$, covers $\conf\times\conf$ (except possibly the diagonal),
and such that
\[
\sigma_t(\{q:(q,q') \in S_t^-\},q') =0 \quad\text{ for
$\measure_t$-almost-every
$q'$.}
\]
Indeed, by the anti-symmetry of $\current_t$,
its positive and negative part $\current_t^+$ and $\current_t^-$
admit supports $S_t^+$ and $S_t^-$ that are transpositions
of each other, whence the result follows.

Put more simply, the mechanism is this: When the current
$\current_t(\D q,\D q')$ is positive, meaning that there should be a
net flow from $\D q'$ to $\D q$, then $\sigma_t(\D q,q') >0$ and
$\sigma_t(\D q',q) =0$, i.e., only jumps from $q'$ to $q$ are allowed;
the converse holds in the case of a negative current. Under all rates
with $j_t = \current_t$, the minimal rates \eqref{tranrates} are
characterized by this property.
\end{enumerate}

\section{Bohmian Mechanics and Bell-Type QFT}

In this section we discuss three particular instances in which
jump rates of the form \eqref{tranrates} play a significant role.

\smallskip
\emph{A.~Bohmian mechanics as continuum limit of jump processes}.
Consider nonrelativistic quantum mechanics for $N$ particles: the
configuration space is $\conf = \RRR^{3N}$, the Hilbert space
$\Hilbert = L^2(\RRR^{3N},\CCC^k)$ and the Hamiltonian
\begin{equation}\label{qmH}
      H = -\sum_{i=1}^N \frac{\hbar^2}{2m_i} \Laplace_i + V(\vr_1,
      \ldots, \vr_N)
\end{equation}
with $\Laplace_i$ the Laplacian acting on the variable $\vr_i$, $m_i$
the mass of the $i$-th particle, and $V$ the potential function
(possibly having values in the Hermitian $k \times k$ matrices).  We
obtain a Markov process on the configuration space in the following
way: first discretize space, i.e., replace $\RRR^3$ by a lattice
$\Lambda = \varepsilon \ZZZ^3$ and the Laplacian $\Laplace_i$ by the
corresponding lattice Laplacian $\Laplace^\varepsilon_i$.  We then can
consider the jump process $Q_t^\varepsilon$ on $\Lambda^N$ with rates
\eqref{tranrates}. As the lattice shrinks, $\varepsilon \to 0$, one
obtains \cite{Sudbery,Vink} in the limit the deterministic
process $Q_t$ satisfying the ordinary differential equation
\begin{equation}\label{Bohm}
      \frac{\D\vQ_{t,i}}{\D t} = \frac{\hbar}{m_i} \im \, \frac{\Psi^*_t
      \nabla_i \Psi_t}{\Psi^*_t \, \Psi_t}(\vQ_{t,1}, \ldots,
      \vQ_{t,N})\,.
\end{equation}
Here $\vQ_{t,i}$ is the $i$-th component of $Q_t$, i.e., the position
of the $i$-th particle, $\Psi_t$ obeys the Schr\"odinger equation
\eqref{schr} with Hamiltonian \eqref{qmH}, and $\Phi_1^* \, \Phi_2$ is
the inner product in $\CCC^k$. The process \eqref{Bohm} is known as
\emph{Bohmian mechanics}
\cite{Bohm52,Bellhidden,survey,DetlefBuch}.  For a suitable other
choice of jump rates \cite{Vink}, also satisfying $j_t = \current_t$
but greater than minimal, one obtains in the continuum limit
$\varepsilon \to 0$ the diffusion process introduced by E.~Nelson and
known as \emph{stochastic mechanics} \cite{stochmech}.

What makes Bohmian mechanics (or, for that matter, stochastic mechanics)
particularly interesting to quantum physicists is that in a Bohmian
universe -- one in which the particles move according to \eqref{Bohm}
and the initial configuration is chosen according to the $|\Psi|^2$
distribution -- the inhabitants find all their observations in
agreement with the probabilistic predictions of quantum mechanics --
in sharp contrast with the traditional belief that it be impossible to
explain the probabilities of quantum mechanics by any theory
describing events objectively taking place in the outside world.

\smallskip
\emph{B.~Bell's jump process for lattice QFT}.  The study of jump
processes with rates \eqref{tranrates} has been inspired by Bohmian
mechanics, in particular by the wish for a theory similar to Bohmian
mechanics covering quantum field theory. The first work in this
direction was Bell's seminal paper \cite{BellBeables}. For
simplicity, Bell replaces physical 3-space by a lattice
$\Lambda$ and considers a QFT on that lattice. A configuration is
specified in his model by the number of fermions $q(\vr)$ at every
lattice site $\vr$. Thus, with the notation $\ZZZ_+ =
\{0,1,2,\ldots\}$, the configuration space is
\[
      \conf = \Gamma(\Lambda) := \Bigl\{ q \in \ZZZ_+^\Lambda:
      \sum_{\vr \in \Lambda} q(\vr) < \infty \Bigr\},
\]
the space of all configurations of a variable (but finite) number of
identical particles on the lattice. (In fact, he imposes a bound on
the total number of particles and assumes that $\Lambda$ is
finite, but this is not necessary.) The PVM $\pov(\,\cdot\,)$ that he
suggests arises from the joint spectral decomposition of the fermion
number operators $N(\vr)$ for every lattice site, i.e., $\pov(q):=
\pov(\{q\})$ is the projection to the joint eigenspace of the
(commuting) operators $N(\vr)$ for the eigenvalues $q(\vr)$. The jump
rate Bell uses is the appropriate special case of \eqref{tranrates}:
the rate of jumping from $q'$ to $q$ is
\begin{equation}\label{Bellrates}
      \sigma_t(q,q') = \frac{\bigl[ \tfrac{2}{\hbar} \im \,
      \sp{\Psi_t}{\pov(q) H \pov(q') \Psi_t} \bigr]^+}
      {\sp{\Psi_t}{\pov(q') \Psi_t}}.
\end{equation}
For studies of Bell's process we refer to
\cite{Sudbery,Vink,BD,Colin,Colinb,crea2,crex1}, and for some
numerical simulations and applications to \cite{Dennis,Dennis2}. We
return to it in more detail in the next section.

\smallskip
\emph{C.~Bohmian mechanics with variable number of particles}.  A
third example of a process for a QFT was considered in \cite{crea1}.
It arose from an attempt to include particle creation and annihilation
into Bohmian mechanics by simply introducing the possibility that
world lines of particles can begin and end. That is, the aim is to
provide a generalization of the Bohmian motion \eqref{Bohm} to a
configuration space of a variable number of particles. Here we
describe this model in a simplified version. For the numerous
similarities between our model process and Bell's discrete process, we
called it a ``Bell-type QFT.''  In \cite{crea2}, methods are developed
for obtaining a canonical Bell-type process for more or less any
regularized QFT.

A configuration of finitely many identical particles can be described
by a finite counting measure on $\RRR^3$.  Since the coincidence
configurations, those in which there are two or more particles at the
same location, form a subset of codimension 3, they are basically
irrelevant, and it will be convenient to exclude them from the
configuration space. What remains, as the space of ``simple
configurations'', is the set of all finite subsets of $\RRR^3$,
\[
      \Gommo (\RRR^3) = \bigl\{ q \subset \RRR^3 : \# q <
      \infty \bigr\}.
\]

Under the physical conditions prevailing in everyday life, the most
frequent type of particle creation and annihilation is the emission
and absorption of photons by electrons. This can be described in a
model QFT as follows.  Particles (photons) move in a Bohmian way and
can be emitted and absorbed by another kind of particles
(electrons). For simplicity, we will assume here that the electrons
remain at fixed locations, given by a finite set $\eta \subset
\RRR^3$; the case of moving electrons is described in \cite{crea1}.
The configuration space is thus the space of photon configurations,
$\conf = \Gommo(\RRR^3)$, and $\Psi_t$ a square-integrable
complex-valued function on $\conf$; the Hilbert space $\Hilbert$ of
these functions is known as the bosonic Fock space arising from
$L^2(\RRR^3)$.

The Markov process $Q_t$ in $\conf$ has piecewise smooth paths. It
obeys the deterministic motion \eqref{Bohm}, interrupted by stochastic
jumps. The process is piecewise deterministic in the sense that,
conditional on the times of two subsequent jumps and the destination
of the first, the path in between these jumps is deterministic. The
jumps correspond to creation or annihilation of a photon near
some point of $\eta$; in particular, every jump changes the number of
photons by one. The process is thus a special kind of a spatial
birth-and-death process with moving individuals \cite{preston}.

The deterministic motion, during which the number of photons is kept
constant, is defined by \eqref{Bohm}; for simplicity, we deviate a
little from the physical facts and assume that the ``photons'' have a
positive mass $m_\mathrm{ph}$. The rate for the transition $q \to q
\cup \vr:= q\cup\{x\}$, i.e., the creation of a new photon at the
location $\vr \in \RRR^3 \setminus q$, has density (with respect to
Lebesgue measure $\D \vr$)
\[
      \sigma_t ( q \cup \vr,q) = \frac{ \Bigl[ \tfrac{2}{\hbar}
      \,\im \, \Psi_t^*(q \cup \vr) \, (\# q +1)^{1/2} \sum\limits_{\vy
      \in \eta} \varphi(\vr-\vy) \, \Psi_t( q) \Bigr]^+}{\Psi^*_t(q) \,
      \Psi_t(q)}
\]
where $\varphi: \RRR^3 \to \RRR$ is a fixed function, a spherically
symmetric, square-integrable potential supported by the ball of radius
$\delta>0$. Likewise, for any $\vr \in q$ the rate for the transition
$q \to q \setminus \vr:=q \setminus \{\vr\}$, i.e., annihilation of
the photon at $\vr$, is
\[
      \sigma_t (q \setminus \vr, q) = \frac{ \Bigl[ \tfrac{2}{\hbar}
      \,\im \, \Psi_t^*(q \setminus \vr) \, (\# q)^{-1/2}
      \sum\limits_{\vy \in \eta} \varphi(\vr-\vy) \, \Psi_t(q)
      \Bigr]^+}{\Psi^*_t(q) \, \Psi_t(q)}
\]
These rates, together with vanishing rate for any other transition,
are in fact a special case of \eqref{tranrates}, for a suitable
integral operator $H_I$ in place of $H$. For the definition of $H_I$
and the derivation of the above expressions from \eqref{tranrates} we
refer to Sect.~3.12 of \cite{crea2}.

Now, $H_I$ is not the Hamiltonian of the relevant QFT, but its
\emph{interaction part}; i.e., the complete Hamiltonian is $H= H_0 +
H_I$, where $H_0$, the \emph{free Hamiltonian}, is given by
\[
      H_0 \Psi(q) = - \sum_{i=1}^{\# q} \frac {\hbar^2}
      {2m_{\mathrm{ph}}} \, \Laplace_{i} \Psi (q).
\]
Observe that there is a correspondence between the splitting $H = H_0
+ H_I$ and the two constituents of the process, the motion \eqref{Bohm}
and the jump rates just given. Deterministic motion corresponds to
$H_0$ while the jumps correspond to $H_I$. Indeed, the \emph{minimal
process}, the one arising as a limiting case from jump processes with
minimal rates \eqref{tranrates}, associated with $H_0$ alone is the
continuous motion \eqref{Bohm} while the minimal process associated
with $H_I$ is the pure jump process with the above rates.

This is an instance of the general rule of \emph{process additivity}:
If the minimal process associated with $H_1$ has generator
$\generator_{1,\Psi_t}$ and the one associated with $H_2$ has
generator $\generator_{2,\Psi_t}$, then the minimal process associated
with $H_1 + H_2$ has generator $\generator_{1,\Psi_t} +
\generator_{2,\Psi_t}$, provided that the (formal) integral kernels of
$H_1$ and $H_2$ have disjoint supports.  The sum of the generators of
a deterministic flow and of a pure jump process generates the
piecewise deterministic process that follows the flow between
stochastic jumps.  In QFT, it is a typical situation that $H = H_0 +
H_I$ where $H_0$ is a differential operator associated with continuous
motion while $H_I$ is an integral operator (often linking different
particle numbers) associated with jumps.

\section{Global Existence of Bell's Jump Process}

In this section we deal with Bell's jump process introduced  as model
$B$
in the last section. As we have shown in \cite{crex1}, this process
exists globally in time. In fact,
for our proof it is not relevant whether the configuration
corresponds to the fermion number operators. We only need that
$\conf$ is any countable set, and $\pov(\,\cdot\,)$ a PVM on $\conf$
acting on $\Hilbert$. In fact we can allow that $\pov(\,\cdot\,)$ is a
positive-operator-valued measure (POVM), a concept somewhat weaker
than a PVM.\footnote{That is, $\pov$ takes values in the
positive (bounded, self-adjoint) operators on $\Hilbert$ (instead of
the projection operators as a PVM) and shares the
countable additivity and normalization of a PVM.}
Here is our result.

\begin{theorem}
      Let $\Hilbert$ be a Hilbert space, $H$ a self-adjoint operator on
      $\Hilbert$, $\conf$ a countable set, and $\pov(\,\cdot\,)$ a POVM 
on
      $\conf$ acting on $\Hilbert$.  For every initial state vector
      $\Psi_0$ with $\|\Psi_0\|=1$ satisfying
      \begin{subequations}\label{assumptions}
      \begin{align}
        \Psi_t \in \domain(H) \quad & \forall\, t \in \RRR,
        \label{assumptiona}\\
        \pov(q) \Psi_t \in \domain(H) \quad & \forall\, t \in \RRR,
        \, q \in \conf,\label{assumptionb}\\
        \int\limits_{t_1}^{t_2} dt \sum_{q,q' \in \conf} \bigl|
        \sp{\Psi_t}{ \pov(q) H \pov(q') \Psi_t} \bigr| < \infty &
        \quad \forall\, t_1,t_2 \in \RRR \text{ with } t_1 < t_2,
        \label{assumptionc}
      \end{align}
      \end{subequations}
      there exists a (right-continuous) Markovian pure jump process
      $(Q_t)_{t \in \RRR}$ on $\conf$ with transition rates
      \eqref{Bellrates} such that, for every $t$, $Q_t$ has distribution
      $\measure_t=\sp{\Psi_t}{ \pov(\,\cdot\,) \Psi_t}$. This process is
unique in
      distribution.
\end{theorem}

Some comments are necessary. First of all, how could the
process fail to exist globally in time? Two kinds of catastrophes could
occur. On the one hand, the jump rate \eqref{Bellrates} is singular
at the nodes of $\Psi$ (i.e., at such $q$ and $t$  for which
$\sp{\Psi_t}{\pov(q) \Psi_t} =0$). While $Q_t$ is sitting on a
configuration $q$ it might become a node, and then the process would not
know how to proceed. It turns out that this
problem does not arise because, with probability one, there is no $t$ at
which $Q_t$ is a node. This is because the
increase of the rates close to the nodes has the positive effect
of forcing the process to jump away
before the singularity time is reached.

The second kind of possible catastrophe would be an explosion,
i.e., an accumulation of infinitely many jumps in finite time.  The
main task is to show that this does not occur, with probability
one. The standard criteria for non-explosion of pure jump processes
are confined to transition rates that are homogeneous in time, relying
heavily on the fact that the holding times are then exponentially
distributed and independent conditionally on the positions; see, e.g.,
Sect.~2.7 of \cite{norris}.  This conditional independence, however,
fails to hold in the case of time-dependent jump rates, and the
singularities of Bell's transition rates do not allow any a priori
bounds as they were used, e.g., in \cite{preston,reuter} to exclude
explosion. The only thing one knows is that the process is designed to
have the prescribed quantum distribution \eqref{measuredef} at fixed
(deterministic) times, and it is this fact we will exploit.  We will
sketch our main arguments below.

Steps towards a proof of global existence of Bell's process have also
been made by G.~Bacciagaluppi \cite{Bthesis,BD}; his approach is,
however, very different from ours.

Concerning the technical assumptions \eqref{assumptions} on $H$,
$\pov$, and $\Psi_0$ we note the following.  For fixed $H$ and $\pov$,
the conditions \eqref{assumptions} define a set of ``good'' initial
state vectors $\Psi_0$ for which we can prove global existence; this
set is obviously invariant under the time evolution \eqref{evolution}.
In fact, when $H$ is a Hilbert--Schmidt operator (i.e., $\tr \, H^2 <
\infty$), the conditions \eqref{assumptions} are satisfied for
\emph{all} POVMs $\pov$ and \emph{all} $\Psi_0 \in \Hilbert$; this is
also true when $H$ is bounded and $\conf$ is finite. (Usually, the
Hamiltonian of a lattice QFT is not Hilbert--Schmidt but can, at
least, be approximated by Hilbert--Schmidt operators. And it is not
unusual in quantum field theory that Hamiltonians need to be ``cut
off'' in one way or another to make them treatable, or well-defined at
all.)  Condition \eqref{assumptionb} ensures that the expression
$\pov(q) H \pov(q') \Psi_t$ can be formed, and thus that
\eqref{Bellrates} is well-defined whenever $q'$ is not a node.

The main construction is obvious. Starting from any fixed initial time
$t_0$, the process $Q_t$ can be expressed for $t>t_0$ in terms of
$T_k$ and $X_k$, the time and the destination of the $k$-th jump after
time $t_0$, via
\[
      Q_t = X_k \quad \text{if } T_k \leq t < T_{k+1}
\]
with $T_0 := t_0$ and $X_0 = Q_{t_0}$. The variables $T_k$ and $X_k$
are defined by their conditional distributions:
\begin{multline}\nonumber
      \prob(T_{k+1}\in\D t, X_{k+1}= q| T_0,X_0,\ldots,T_k,X_k) =\\
      \1_{\{T_k < t\}} \, \sigma_t(q, X_k) \, \exp \Bigl(
      -\int\limits_{T_k}^t  \sigma_s(\conf,X_k)\,\D s \Bigr) \D t,
\end{multline}
where the role of the ``failure rate function'' is played by
\[
      \sigma_t(\conf,q) = \sum_{q' \in \conf} \sigma_t(q',q)\,,
\]
the total jump rate, to whatever destination, at $q$.

Here is the reason why the process cannot run into a node.  By
definition, the conditional probability of remaining at $q$ until at
least time $t_2$, given that $Q_{t_1} = q$, is
\[
      \exp\Bigl(-\int\limits_{t_1}^{t_2}\sigma_t(\conf,q)  \, \D t 
\Bigr).
\]
We want to show that this probability vanishes whenever $q$ is
a node at $t_2$ but not at any $t$ with $t_1 \leq t < t_2$. Ignoring
some technical subtleties, this can be derived by the following simple
calculation. Since a sum of positive parts exceeds the positive part
of the sum, we conclude from \eqref{Bellrates} that
\[
      \sigma_t(\conf,q') =  \sum\limits_{q\in\conf}
      \frac{[\tfrac{2}{\hbar} \, \im \, \sp{\Psi_t}{\pov(q) H \pov(q')
      \Psi_t}]^+}{\sp{\Psi_t}{\pov(q') \Psi_t}} \geq
      \frac{[\tfrac{2}{\hbar} \, \im \, \sp{\Psi_t}{H \pov(q')
      \Psi_t}]^+}{\sp{\Psi_t}{\pov(q') \Psi_t}}\,.
\]
Omitting the positive part and using \eqref{dotmeasure} we find
\[
      \sigma_t(\conf,q')  \geq  \frac{-(\D/\D t)
      \sp{\Psi_t}{ \pov(q') \Psi_t}} {\sp{\Psi_t}{\pov(q') \Psi_t}} = -
      \frac{\D}{\D t} \log \sp{\Psi_t}{\pov(q') \Psi_t}
\]
at every $t$ with $\sp{\Psi_t}{\pov(q') \Psi_t}>0$. Hence, by the
fundamental theorem of calculus,
\[
      \int\limits_{t_1}^{t_2} \sigma_t(\conf,q) \, \D t \geq - \log
      \sp{\Psi_{t_2}} {\pov(q) \Psi_{t_2}} + \log \sp{\Psi_{t_1}}
      {\pov(q) \Psi_{t_1}} = \infty,
\]
since $q$ is a node at $t_2$ (so that the first term is infinite) but
not at $t_1$ (so that the second term is finite).

Another part of the proof we would like to sketch here is the core of
the argument why the jump times cannot accumulate. As a convenient
notation, we introduce an additional ``cemetery'' configuration
$\friedhof$ and set $Q_t := \friedhof$ for all $t$ after the explosion
time $\sup_k T_k$. Let $S(t_1,t_2)$ be the number of jumps that the
process performs in the time interval $[t_1,t_2]$.  The random
variable $S(t_1,t_2)$ is either a nonnegative integer or infinite. Our
assumption \eqref{assumptionc} implies in fact that $S(t_1,t_2)$ has
finite expectation, for all $t_1<t_2$, and thus is finite almost
surely. To see this we use the equation
\begin{equation}\label{Ssigmarho}
      \EEE\, S(t_1,t_2) = \int\limits_{t_1}^{t_2}  \sum_{q,q' \in \conf}
      \sigma_t(q,q') \, \rho_t(q')\,\D t,
\end{equation}
with $\rho_t(q') = \prob(Q_t = q')$. Intuitively, equation
\eqref{Ssigmarho} can be understood as follows. $\sigma_t(q,q') \,
\rho_t(q') \, \D t$ is the probability of a jump from $q'$ to $q$ in
the infinitesimal time interval $[t,t+\D t]$. Summing over $q$ and
$q'$, we obtain (the expected total jump rate and thus) the total
probability of a jump during $[t,t+\D t]$. Integrating over $t$, we
obtain the expected number of jumps. The point of equation
\eqref{Ssigmarho} is that its right-hand side involves exclusively the
one-time quantities $\sigma_t$ and $\rho_t$. Now, one can deduce from
the definition of the process that
\begin{equation}\label{rholeqPsi2}
      \rho_t(q) \leq \sp{\Psi_t}{\pov(q) \Psi_t},
\end{equation}
where any case of strict inequality would have to go along with a
positive probability $\prob(Q_t = \friedhof)$ of accumulation before
$t$.
Combining \eqref{Ssigmarho} and \eqref{rholeqPsi2}, we obtain
\[
\begin{split}
      \EEE\, S(t_1,t_2) &\leq \int\limits_{t_1}^{t_2}\sum_{q,q' \in
      \conf} \sigma_t (q,q') \, \sp{\Psi_t}{\pov(q') \Psi_t}\, \D t \\
       &= \int\limits_{t_1}^{t_2} \sum_{q,q' \in \conf} \bigl[
      \tfrac{2}{\hbar} \, \im \,\sp{\Psi_t}{\pov(q) H \pov(q') \Psi_t}
      \bigr]^+  \,\D t\\
       &\leq \tfrac{2}{\hbar} \int\limits_{t_1}^{t_2} \sum_{q,q' \in
      \conf} \bigl| \sp{\Psi_t}{\pov(q) H \pov(q') \Psi_t} \bigr| \, \D t
<
\infty
      \end{split}
\]
by assumption \eqref{assumptionc} of Theorem~1. Indeed, this
reasoning more or less dictates the assumption \eqref{assumptionc}.

\section{Other Global Existence Questions}

Variants of the reasoning in the previous paragraph could be used in
the future in other global existence proofs. Here is an example
concerning the \emph{global existence of Bohmian mechanics} (with
fixed number of particles). This was first proved in \cite{bmex} under
suitable assumptions on the potential $V$ and the initial wavefunction
$\Psi_0$. One way in which a solution $Q_t$ of \eqref{Bohm} may fail
to exist globally is by escaping to infinity (i.e., leaving every
bounded set in $\RRR^{3N}$) in finite time. To control this behavior,
we suggest to consider an analogue of $S(t_1,t_2)$: Let
$D(t_1,t_2)$ be the Euclidean distance in $\RRR^{3N}$ traveled by
$Q_t$ between $t_1$ and $t_2$. Then
\begin{equation}\label{Dvrho}
      \EEE\, D(t_1,t_2) = \int\limits_{t_1}^{t_2} \D t
      \int\limits_{\RRR^{3N}} \D q \, |v_t(q)| \, \rho_t(q),
\end{equation}
where $v_t$ is the Bohmian velocity vector field on $\RRR^{3N}$, with
$i$-th component given by the right-hand side of \eqref{Bohm}, and the
expectation is taken over the randomness coming from the initial
configuration.  If this expectation can be shown to be finite,
$D(t_1,t_2)$ must be almost surely finite.  Using $\rho_t(q) \leq
|\Psi_t(q)|^2$, the analogue of \eqref{rholeqPsi2}, we see that the
escape to infinity is almost surely excluded provided that
\[
      \int\limits_{t_1}^{t_2} \D t \int\limits_{\RRR^{3N}} \D q \,
      |\Psi_t^*(q) \, \nabla \Psi_t(q)| < \infty \quad \text{whenever }
      t_1< t_2.
\]
This is a condition analogous to \eqref{assumptionc}; it is almost
equivalent to the assumption A4 of \cite{bmex}. The proof there,
however, is different, using skillful
estimates of the probability flux across suitable surfaces in
configuration space-time $\RRR^{3N} \times \RRR$ surrounding the
``bad'' points (nodes, infinity, points where $\Psi$ is not
differentiable). The above argument based on \eqref{Dvrho} might
contribute to a simpler global existence proof.

The global existence question is still open for the \emph{Bohm--Dirac
law of motion} \cite{Bohm53,BH}, a version of Bohmian mechanics
suitable for wavefunctions $\psi$ obeying the Dirac equation. The
Dirac equation is a relativistic version of the Schr\"odinger equation
and reads
\begin{equation}\label{Dirac}
      \I\hbar\frac{\partial \psi}{\partial t} = -\sum_{i=1}^N \bigl(\I c
      \hbar \valpha_i \cdot \nabla_i \psi + \beta_i m_i c^2 \psi \bigr)
\end{equation}
where $\psi_t$ is a function $\RRR^{3N} \to (\CCC^{4})^{\otimes N}$,
$c$ is the speed of light, $m_i$ the mass of the $i$-th particle,
$\valpha_i$ the vector of Dirac alpha matrices acting on the $i$-th
spin index of the wavefunction, and $\beta_i$ the Dirac beta matrix
acting on the $i$-th index.  The Bohm--Dirac equation of motion reads
\begin{equation}\label{BohmDirac}
      \frac{\D\vQ_{t,i}}{\D t} = c\frac{\psi^*_t \valpha_i
      \psi_t}{\psi^*_t \psi_t} (\vQ_{t,1}, \ldots, \vQ_{t,N})
\end{equation}
where $\phi^*\psi$ is the inner product in Dirac spin-space. Since
these velocities are bounded by the speed of light, the question of
escape to infinity does not arise here. But the question of running
into a node does because, like the minimal jump rate
\eqref{tranrates}, the velocity formula \eqref{BohmDirac} is
ill-defined at nodes.

This question can be treated in a way analogous to the previous
arguments based on \eqref{Ssigmarho} and \eqref{Dvrho}. To this end,
let $L(t_1,t_2)$ be the total variation, between $t_1$ and $t_2$, of
$\log |\psi_t(Q_t)|^2$; in other words,
\[
      L(t_1,t_2) = \int\limits_{t_1}^{t_2} \D t \, \Bigl| \frac{\D}{\D t}
      \log |\psi_t(Q_t)|^2 \Bigr|.
\]
It takes values in $[0,\infty]$ and is infinite if the trajectory
$Q_t$ runs into a node between $t_1$ and $t_2$. This must be a null
event if $L(t_1,t_2)$ has finite expectation; for the latter we have
the formula
\[
      \EEE\, L(t_1,t_2) = \int\limits_{t_1}^{t_2} \D t
      \int\limits_{\RRR^{3N}} \D q \, \Bigl| \Bigl( \frac {\partial}
      {\partial t} + v_t(q) \cdot \nabla\Bigr) \log |\psi_t(q)|^2 \Bigr|
      \, \rho_t(q),
\]
where $v_t$ is the vector field on $\RRR^{3N}$ whose $i$-th component
is the right-hand side of \eqref{BohmDirac}.  Using $\rho_t(q) \leq
|\psi_t(q)|^2$ we obtain
\[
      \EEE\, L(t_1,t_2) \leq \int\limits_{t_1}^{t_2} \D t
      \int\limits_{\RRR^{3N}} \D q \, \Bigl| \Bigl( \frac {\partial}
      {\partial t} + v_t(q) \cdot \nabla\Bigr) |\psi_t(q)|^2 \Bigr| \,.
\]
As the velocities are bounded, the last expression is at most
\[
     \int\limits_{t_1}^{t_2} \D t \int\limits_{\RRR^{3N}} \D q \,
      \biggl( \Bigl| \frac {\partial} {\partial t} |\psi_t(q)|^2 \Bigr| +
      c\Bigl| \nabla |\psi_t(q)|^2 \Bigr| \biggr)\,.
\]
Inserting \eqref{Dirac} and using that the alpha and beta matrices
have norm 1 we see that this in turn is not larger than
\[
     \int\limits_{t_1}^{t_2} \D t \int\limits_{\RRR^{3N}} \D q \,
      \biggl( 2c\Bigl| \psi^*_t \nabla \psi_t \Bigr| + \tfrac{2}{\hbar}
      \Bigl(\sum_i m_i c^2\Bigr) \psi_t^* \psi_t + 2c\Bigl| \psi_t^*
      \nabla \psi_t \Bigr| \biggr)\,.
\]
Since $\|\psi_t\|=1$, the last integral coincides with
\[
      \tfrac{2}{\hbar} \Bigl(\sum_i m_i c^2\Bigr) (t_2-t_1) + 4c
      \int\limits_{t_1}^{t_2} \D t \int\limits_{\RRR^{3N}} \D q \, \Bigl|
      \psi^*_t \nabla \psi_t \Bigr|\,.
\]
The question remains under which conditions on $\psi_0$ the last term
is finite for arbitrary $t_1 < t_2$. This is presumably the case when
$\psi_0$ lies in Schwartz space (containing all smooth functions $f$
such that $f$ and all its derivatives decay, at infinity, faster than
$|q|^{-n}$ for any $n>0$). To work out the proof remains for future
research \cite{TT04}.

While global existence of the Bohm--Dirac trajectories can presumably
be proved with the same methods as used in \cite{bmex} for
\eqref{Bohm} (estimating the flux across suitable surfaces surrounding
the bad points), such a proof requires a lot of effort. It seems that
the argument just sketched would be much easier and more elegant.

Another global existence question that is still open is that
concerning the process defined in \cite{crea1} and described above in
Part~C of Sect.~3, and for similar processes, on a configuration space
that is a disjoint union of manifolds, following deterministic
trajectories interrupted by stochastic jumps.

\section{Deterministic Jumps and Boundaries in Configuration Space}

In this last section we describe another application of minimal jump
rates, one that has not yet been discussed in the literature and that
raises some questions for further research. Suppose that the
configuration space $\conf$ is a Riemannian manifold with boundaries,
or more generally the disjoint union of (at most countably many)
Riemannian manifolds with boundaries. We write $\conf = \boundary \cup
\interior$ where $\boundary$ denotes the boundary and $\interior$ the
interior.

We develop below an analogue of Bohmian mechanics on $\conf$,
consisting of smooth motion interrupted by jumps from the boundary to
the interior or vice versa.  The jumps from $\boundary$ to $\interior$
are deterministic and occur whenever the process hits the
boundary. The jumps from $\interior$ to $\boundary$ are stochastic,
and their rates are fully determined by requiring that (i) the process
is Markovian and equivariant, and (ii) the construction is invariant
under time reversal, in that the processes associated with $\Psi_t$
resp.\ $\Psi_{-t}^*$ are reverse to each other, in distribution. These
rates are, in fact, another instance of minimal jump rates.

Configuration spaces with boundaries arise from QFT if a particular
kind of ``ultraviolet cutoff'' is applied, which can be regarded as
corresponding to smearing out the charge of an electron over a sphere
rather than a ball. Here is an example.  Consider again electrons and
photons, with the electrons fixed at locations given by the finite set
$\eta \subset \RRR^3$, and suppose that photons cannot get closer than
a fixed (small) distance $\delta>0$ to an electron, as they get
absorbed when they reach that distance.  Thus, the available
configuration space is
\begin{equation}\label{confdef}
      \conf = \bigl\{ q \in \Gommo(\RRR^3) : d(q,\eta)\geq \delta
      \bigr\}\,,
\end{equation}
where $d(q,\eta)$ is the Euclidean distance of the finite sets $q$ and
$\eta$.  The space $\conf$ is a countable disjoint union of Riemannian
manifolds with boundary,
\[
      \conf= \bigcup_{n=0}^\infty \bigl\{ q \in \conf : \#q=n \bigr\}.
\]
Its interior is $\interior = \{ q \in \Gommo(\RRR^3) : d(q,\eta)>
\delta \}$, and its boundary $\boundary =\{ q \in
\Gommo(\RRR^3) : d(q,\eta)= \delta \}$.

For this or any other configuration space with boundaries,
the law of motion
\begin{equation}\label{Bohm2}
      \frac{\D Q_t}{\D t} = v_t(Q_t) = \frac{\hbar}{m} \, \im \,
      \frac{\Psi_t^* \nabla \Psi_t}{\Psi_t^* \Psi_t}(Q_t)
\end{equation}
on $\interior$ must be completed by specifying what should happen to
the process at the time $\tau$ when it reaches the boundary. (No
specification is needed, however, for what should happen when two
photons collide, as this has probability zero ever to occur.)  The
specification we consider here is a deterministic jump law
\[
      Q_{\tau+} = \dest(Q_{\tau-})
\]
for a fixed mapping $\dest:\boundary \to \interior$.  In our
example \eqref{confdef}, the obvious choice of $\dest$ is
\[
      \dest(q) = \{\vr \in q: d(\vr,\eta)> \delta  \},
\]
which means that all photons having reached the critical distance
$\delta$ to some electron disappear.

Since we want the theory to be reversible, we must also allow for
spontaneous jumps from interior points to boundary points.  Since we
want the process to be an equivariant Markov process, the rate for a
jump from $q' \in \interior$ to a surface element $\D q \subseteq
\boundary$ must be, as one can derive,
\begin{equation}\label{boundrate}
      \sigma_t(\D q,q') = \frac{\bigl[ n(q) \cdot v_t(q) \, |\Psi_t(q)|^2
      \bigr]^+} {|\Psi_t(q')|^2} \, \antidest(\D q,q'),
\end{equation}
where $n(q)$ is the inward unit normal vector to the boundary at $q
\in \boundary$, the dot $\cdot$ denotes the Riemannian inner product,
and $\antidest(\D q,q')$ is the measure-valued function defined in
terms of the Riemannian volume measure $\vol$ on $\conf$ and the
Riemannian surface area measure $\area$ on $\boundary$ by
\[
    \antidest(B,q') = \frac{\area(B \cap \dest^{-1}(\D q'))}{\vol(\D 
q')},
\]
with $B \subseteq \boundary$, and the right-hand side denoting a
Radon--Nikodym derivative (the existence of which we presuppose).  The
measure $\antidest(\,\cdot\,,q')$ is concentrated on the subset
$\dest^{-1}(q')$ of the boundary for almost every $q'$.  For a
probability distribution on $\conf$ having a density function $\rho$
with respect to $\vol$, one obtains the following
probability transport equation at $q' \in \interior$:
\begin{equation}\label{continuity2}
      \frac{\partial \rho_t}{\partial t}(q') = -\nabla \cdot \bigl( 
\rho_t
      v_t \bigr) (q') -\sigma_t(\boundary,q') \, \rho_t(q') + \!
      \int\limits_{\boundary} \!\antidest(\D q,q') \Bigl[n(q) \cdot 
v_t(q)
      \, \rho_t(q) \Bigr]^- \!\! .
\end{equation}
For equivariance we need that \eqref{continuity2}, with $|\Psi_t|^2$
in place of $\rho_t$, has the structure of the transport equation for
$|\Psi_t|^2$ that follows from the Schr\"odinger equation
\eqref{schr},
\begin{equation}\label{dPsi2dt}
      \frac{\partial |\Psi_t(q')|^2}{\partial t} = \tfrac{2}{\hbar} \, 
\im
      \, \Psi^*_t(q') \, (H\Psi_t)(q').
\end{equation}
This is not automatically the case, but it follows from (and
therefore suggests) the following boundary condition relating
$\Psi_t|_{\boundary}$ to $\Psi_t|_{\interior}$: for all $q\in
\boundary$,
\begin{equation}\label{bouncon}
      n(q) \cdot \nabla \Psi_t(q) = \gamma(q) \, \Psi_t(f(q))
\end{equation}
where $\gamma(q)$ is any complex number.\footnote{More generally, if
$\Psi$ takes values not in $\CCC$ but in a higherdimensional complex
vector space, $\gamma(q)$ would be a $\CCC$-linear mapping from the
value space at $f(q)$ to the value space at $q$.}  This condition
prescribes the normal derivative of the wavefunction on the
boundary. Some boundary condition would be needed anyway to define the
evolution of the wavefunction, i.e., to select a self-adjoint
extension of the Laplacian; whether \eqref{bouncon} actually suffices
for this, we have to leave to future research. Note that
\eqref{bouncon} is a linear condition and thus defines a subspace of
the Hilbert space $L^2(\conf)$.  {}From \eqref{bouncon} and
\eqref{Bohm2}, one obtains equivariance with respect to the formal
Hamiltonian $H = -\tfrac{\hbar^2}{2m} \Laplace + H_I$, where
\[
\begin{split}
      \sp{\Phi}{H_I \Psi} =
      & \tfrac{\hbar^2}{2m} \int\limits_{ \conf} \vol(\D q')
      \int\limits_{\boundary} \antidest(\D q,q')\, \Phi^*(q')\,
      \gamma^*(q) \, \Psi(q)
      \\
      &+ \tfrac{\hbar^2}{2m} \int\limits_{\conf} \vol(\D q)
      \int\limits_{\boundary} \antidest(\D q',q)\, \Phi^*(q')\,
      \gamma(q') \, \Psi(q).
\end{split}
\]
Furthermore, the jump rate \eqref{boundrate} is
indeed the minimal jump rate \eqref{tranrates} associated with $H_I$,
thanks to \eqref{bouncon}.

Let us emphasize the following aspects. First, starting from the
picture of a piecewise deterministic process that jumps whenever it
hits the boundary, we arrived with remarkable ease at the probability
transport equation \eqref{continuity2} and thus at the boundary
condition \eqref{bouncon}. Second, we have \emph{derived} what the
interaction Hamiltonian $H_I$ is; once the destination mapping $f$ and
the corresponding coefficients $\gamma$ had been selected, there was
no further freedom. Third, it turned out that the minimal jump rate is
the only \emph{possible} jump rate in this setting. Its very
minimality plays a crucial role: a jump to a boundary point $q$ at
which the velocity field is pointing \emph{towards} the boundary,
$n(q) \cdot v_t(q) < 0$, would not allow any continuation of the
process since there is no trajectory starting from $q$. The problem is
absent if the velocity at $q$ is pointing \emph{away} from the
boundary, $n(q) \cdot v_t(q) >0$. (We are leaving out the case $n(q)
\cdot v_t(q) = 0$.) On the other hand, jumps from $q$ to $f(q)$ cannot
occur when $v_t(q)$ is pointing away from the boundary since in that
case there is no trajectory arriving at $q$. Thus, the jumps must be
such that at each time $t$, one of the transitions $q \to f(q)$ or
$f(q) \to q$ is forbidden, and the decision is made by the sign of
$n(q) \cdot v_t(q)$.


\begin{thebibliography}{[BDDGZ95]}

\bibitem[Bac96]{Bthesis} Bacciagaluppi, G.: Topics in the Modal
      Interpretation of Quantum Mechanics. Ph.~D. thesis, University of
      Cambridge (1996)

\bibitem[BD99]{BD} Bacciagaluppi, G., Dickson, M.: Dynamics for modal
      interpretations. Found. Phys., \textbf{29},
      1165--1201 (1999)

\bibitem[Bel66]{Bellhidden} Bell, J.S.: On the problem of hidden
      variables in quantum mechanics. Rev. Mod. Phys., \textbf{38},
      447--452 (1966). Reprinted in: Bell, J.S.: Speakable and 
unspeakable
      in quantum mechanics.  Cambridge University Press, Cambridge 
(1987),
      p.~1.

\bibitem[Bel86]{BellBeables} Bell, J.S.: Beables for quantum field
      theory. Phys. Rep., \textbf{137}, 49--54 (1986).  Reprinted in:
      Bell, J.S.: Speakable and unspeakable in quantum mechanics.
      Cambridge University Press, Cambridge (1987), p.~173.  Also
      reprinted in: Peat, F.D., Hiley, B.J. (eds) Quantum Implications:
      Essays in Honour of David Bohm. Routledge, London (1987), p.~227.
      Also reprinted in: Bell, M., Gottfried, K., Veltman, M. (eds) John
      S.\ Bell on the Foundations of Quantum Mechanics. World Scientific
      Publishing (2001), chap.~17.

\bibitem[BDGPZ95]{bmex} Berndl, K., D\"urr, D., Goldstein, S.,
      Peruzzi, G., Zangh{\`\i}, N.: On the global existence of Bohmian
      mechanics.  Commun. Math. Phys., \textbf{173}, 647--673 (1995), and
      quant-ph/9503013

\bibitem[BDDGZ95]{survey} Berndl, K., Daumer, M., D\"urr, D.,
      Goldstein, S., Zangh\`\i, N.: A Survey of Bohmian Mechanics. Il
      Nuovo Cimento B, \textbf{110}, 737--750 (1995), and 
quant-ph/9504010

\bibitem[Boh52]{Bohm52} Bohm, D.: A Suggested Interpretation of the
      Quantum Theory in Terms of ``Hidden'' Variables, I and
      II. Phys. Rev., \textbf{85}, 166--193 (1952)

\bibitem[Boh53]{Bohm53} Bohm, D.: Comments on an Article of Takabayasi
      concerning the Formulation of Quantum Mechanics with Classical
      Pictures. Progr. Theoret. Phys., \textbf{9}, 273--287 (1953)

\bibitem[BH93]{BH} Bohm, D., Hiley, B.J.: The Undivided Universe: An
      Ontological Interpretation of Quantum Theory.  Routledge, London
      (1993)

\bibitem[Col03a]{Colin} Colin, S.: The continuum limit of the Bell
      model.  quant-ph/0301119

\bibitem[Col03b]{Colinb} Colin, S.: A deterministic Bell model.  Phys.
      Lett. A, \textbf{317}, 349--358 (2003), and quant-ph/0310055

\bibitem[DR03]{Dennis} Dennis, E., Rabitz, H.: Bell trajectories for
      revealing quantum control mechanisms. Phys. Rev. A, \textbf{67},
      033401 (2003), and quant-ph/0208109

\bibitem[Den03]{Dennis2} Dennis, E.: Purifying Quantum States: Quantum
      and Classical Algorithms. Ph.D. thesis, University of California,
      Santa Barbara (2003)

\bibitem[D\"ur01]{DetlefBuch} D\"urr, D.: Bohmsche Mechanik als
      Grundlage der Quantenmechanik. Springer, Berlin (2001)

\bibitem[DGTZ03a]{crea1} D\"urr, D., Goldstein, S., Tumulka, R.,
      Zangh{\`\i}, N.: Trajectories and Particle Creation and 
Annihilation
      in Quantum Field Theory. J. Phys. A: Math. Gen., \textbf{36},
      4143--4149 (2003), and quant-ph/0208072

\bibitem[DGTZ03b]{crlet} D\"urr, D., Goldstein, S., Tumulka, R.,
      Zangh{\`\i}, N.: Bohmian Mechanics and Quantum Field
      Theory. quant-ph/0303156

\bibitem[DGTZ03c]{crea2} D\"urr, D., Goldstein, S., Tumulka, R.,
      Zangh{\`\i}, N.: Quantum Hamiltonians and Stochastic
      Jumps. quant-ph/0303056

\bibitem[GT03]{crex1} Georgii, H.-O., Tumulka, R.: Global Existence of
      Bell's Time-Inhomogeneous Jump Process for Lattice Quantum Field
      Theory. math.PR/0312294 and mp\underline{\ }arc~04-11

\bibitem[Nel85]{stochmech} Nelson, E.: Quantum Fluctuations.
      Princeton University Press, Princeton (1985)

\bibitem[Nor97]{norris} Norris, J. R.: Markov chains.  Cambridge
      University Press, Cambridge (1997)

\bibitem[Pre76]{preston} Preston, C. J.: Spatial birth-and-death
      processes. Bull. Inst. Internat. Statist., \textbf{46}, no. 2,
      371--391, 405--408 (1975)

\bibitem[RL53]{reuter} Reuter, G. E. H., Ledermann, W.: On the
      differential equations for the transition probabilities of Markov
      processes with enumerably many states. Proc. Cambridge Philos. 
Soc.,
      \textbf{49}, 247--262 (1953)

\bibitem[RS90]{Roy} Roy, S.M., Singh, V.: Generalized beable quantum
      field theory. Phys. Lett. B, \textbf{234}, 117--120 (1990)

\bibitem[Sud87]{Sudbery} Sudbery, A.: Objective interpretations of
      quantum mechanics and the possibility of a deterministic
      limit. J. Phys.  A: Math. Gen., \textbf{20}, 1743--1750 (1987)

\bibitem[TT04]{TT04} Teufel, S., Tumulka, R.: Global Existence of
      Bohm--Dirac Trajectories for Nonsingular Potentials. In
      preparation.

\bibitem[Vin93]{Vink} Vink, J.C.: Quantum mechanics in terms of
      discrete beables. Phys. Rev. A, \textbf{48}, 1808--1818 (1993)

\end{thebibliography}
\end{document}